\documentclass[12pt]{amsart}
\usepackage{graphicx}
\usepackage{graphics}
\usepackage{amsmath}
\usepackage{amscd}
\usepackage{latexsym}

\begin{document}

\textwidth 6.2in
\textheight 7.6in
\evensidemargin .75in
\oddsidemargin.75in

\newtheorem{Thm}{Theorem}
\newtheorem{Lem}[Thm]{Lemma}
\newtheorem{Cor}[Thm]{Corollary}
\newtheorem{Prop}[Thm]{Proposition}
\newtheorem{Rm}{Remark}

\def\a{{\mathbb a}}
\def\C{{\mathbb C}}
\def\A{{\mathbb A}}
\def\B{{\mathbb B}}
\def\D{{\mathbb D}}
\def\E{{\mathbb E}}
\def\R{{\mathbb R}}
\def\P{{\mathbb P}}
\def\S{{\mathbb S}}
\def\Z{{\mathbb Z}}
\def\O{{\mathbb O}}
\def\H{{\mathbb H}}
\def\V{{\mathbb V}}
\def\Q{{\mathbb Q}}
\def\Cn{${\mathcal C}_n$}
\def\CM{\mathcal M}
\def\CG{\mathcal G}
\def\CH{\mathcal H}
\def\CT{\mathcal T}
\def\CF{\mathcal F}
\def\CA{\mathcal A}
\def\CB{\mathcal B}
\def\CD{\mathcal D}
\def\CP{\mathcal P}
\def\CS{\mathcal S}
\def\CZ{\mathcal Z}
\def\CE{\mathcal E}
\def\CL{\mathcal L}
\def\CV{\mathcal V}
\def\CW{\mathcal W}
\def\IC{\mathbb C}
\def\IF{\mathbb F}
\def\IK{\mathcal K}
\def\IL{\mathcal L}
\def\IP{\bf P}
\def\IR{\mathbb R}
\def\IZ{\mathbb Z}

\title{Knot surgery and Scharlemann manifolds}
\author{Selman Akbulut}
\thanks{The author is partially supported by NSF grant DMS 0905917}
\address{Department  of Mathematics, Michigan State University,  MI, 48824}
\email{akbulut@math.msu.edu }
\date{\today}
\begin{abstract} 
We discuss the relation between Fintushel-Stern knot surgery operation on $4$-manifolds and ``Scharlemann manifolds'',  and as a corollary show that they all are standard. Along the way we show the fishtail can exotically knot in $S^4$ infinitely many ways.
\end{abstract}

\date{}
\maketitle

\setcounter{section}{-1}

\vspace{-.1in}

\section{Introduction}

Let X be a smooth 4-manifold, and $T^{2}\times D^2\subset X$ be an imbedded torus with trivial normal bundle, and $K\subset S^3 $ be a knot, $N(K)$ be its tubular neighborhood. The  Fintushel-Stern {\it knot surgery operation}
 is the operation of replacing $T^2 \times D^2$  with $(S^{3} - N (K)) \times S^1$, so that the meridian $p \times \partial D^2$  of the torus coincides with the longitude of $K$ \cite{fs}. 
 $$X \leadsto X_{K}= (X- T^2 \times D^2)\cup (S^3-N(K))\times S^1$$ 

The handlebody picture of this operation was given in \cite{a1}.  Let  $K\subset S^3$ be a knot, and $S^3_{K}$ be the $3$-manifold obtained from $S^3$ by $\pm 1$ surgery  to $K$ (either one).
 The (generalized)  {\it  Scharlemann manifold}  $M(K)$ is the manifold obtained by surgering the circle $C\subset S^{1}\times S^{3}_{K}$ (with even framing) which corresponds to the meridian of the knot $K$. It is clear that $M(K)$ is homotopy equivalent to $S^{1}\times S^{3} \# \;(S^{2} \times S^2)$. 
 in \cite{s} Scharlemann had posed the question whether $M(K)$ is standard when $K$ is the trefoil knot;  and in \cite{a2} this question was answered affirmatively. Here we show that $M(K)$ is also standard for any $K$. 
    We decided to write this paper after seeing \cite{t} which claims the same result.  We felt that there should be a natural direct proof generalizing the steps of \cite{a2} by using the knot-surgery description of  $4$-manifolds \cite{a1}.   It turns out that the stabilization theorem of  \cite{a3} provides the necessary tool to link these two. Along the way we  relate  Scharlemann manifolds $M(K)$ to the knot surgery operation $X \leadsto X_{K}$, and give a sufficient criterion when a knot surgery operation doesn't change the smooth structure of the underlying manifold.
      I thank M.Tange for stimulating my interest to revisit this problem.
  
 \vspace{.05in}
  
\section{A review of the stabilization}

In \cite{a3} (and also in \cite{au}) it was shown that $X_{K}$ is stably trivial, i.e.:
$$X_{K} \# (S^{2}\times S^2 )=X \# (S^2 \times S^2)$$
In  \cite{a3} a specific trivialization move was described in terms of  handles (i.e. turning a ``ribbon $1$-handle'' to $2$-handle'). More specifically  it was shown that surgering the circle $A\subset T^2\times D^2$ (as shown in Figure 1) gives the same manifold as surgering the corresponding $A\subset (T^2\times D^2)_{K}$.

 \begin{figure}[ht]  \begin{center}  
\includegraphics[width=.7\textwidth]{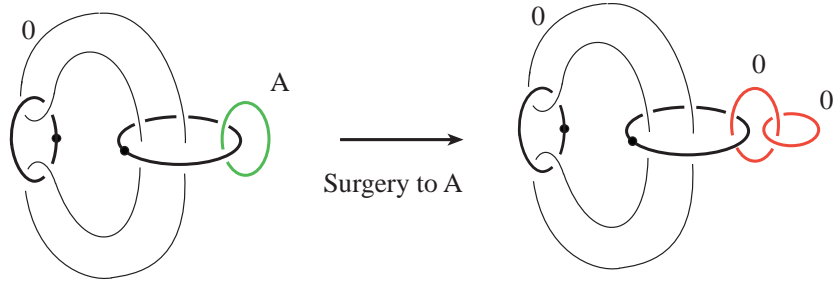}   
\caption{Surgering $T^2\times D^2$} 
\end{center}
\end{figure}

 Notice, if we attach a $2$-handle $h^2$ to $T^2 \times D^2$ along $A$ with zero framing (as in Figure 2),  we get  $ \Gamma:=T^{2}\times D^2  + h^2 =S^{1} \times B^3 \; \natural \; (S^2\times B^2)$, and this identification takes the loop $B$ to the meridian of $S^2 \times B^2$. 
 
  \begin{figure}[ht]  \begin{center}  
\includegraphics[width=.6\textwidth]{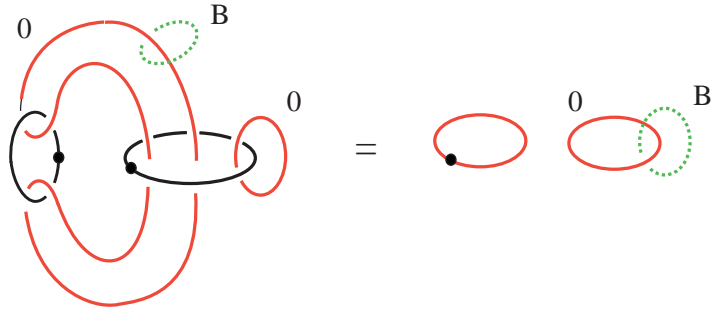}   
\caption{$S^{1} \times B^3 \; \natural \; (S^2\times B^2)$} 
\end{center}
\end{figure}

  The proof of \cite{a3} shows that  the knot surgery of  $S^{1} \times B^3 \; \natural \; (S^2\times B^2)$ along this $T^2\times D^2 \subset S^{1} \times B^3 \; \natural \; (S^2\times B^2)$ keeps it standard:

{\Thm  [\cite{a3}] $[\; S^{1} \times B^3 \; \natural \; (S^2\times B^2)\; ]_{K} = S^{1} \times B^3 \; \natural \; (S^2\times B^2)$ } 

\proof (Sketch) Figure 4 gives the handlebody of the knot surgery (where K is drawn as the trefoil knot). The zero framed linking circle to of  the ``ribbon $1$-handle''  cancels this ribbon $1$-handle,  and in the process the rest of the handlebody becomes standard (cf. \cite{a3}).\qed

\vspace{.05in}

This theorem gives a sufficient condition for showing that a knot surgery operation does not change the underlying smooth manifold.  More specifically,  If  a  torus $T^2\subset X$  has a $\Gamma = S^{1}\times B^3 \#(S^2 \times B^2)$ neighborhood in $X$ (put another way, if the loop $A\subset \partial (T^2\times D^2)$ bounds a disk in the complement $X - T^2\times D^2$, whose normal framing induces the zero framing on  $A$), then $X_{K}=X$. For example $S^4$ can be decomposed as a union of two fishtails glued along boundaries as in Figure 3, and clearly  the torus inside has a $\Gamma $ neighborhood.

 \begin{figure}[ht]  \begin{center}  
\includegraphics[width=.55\textwidth]{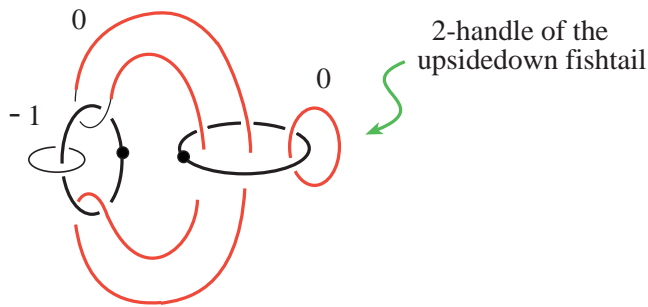}   
\caption{$S^4$ as a union of two fishtails} 
\end{center}
\end{figure} 

Similarly Figure 5  describes of  $S^2\times S^2$ as the double of the cusp, and Figure 6 describes $S^{1}\times S^3 \# (S^2 \times S^2)$  as the double of the fishtail.  Note that  in these figures we give some alternative pictures of these handlebodies by using the diffeomorphism $\varphi: S^{2}\times T^2 \to  S^{2}\times T^2$ of Figure 7, which carries the loop $A$ to itself by twisting its tubular neighborhood. Clearly (from the pictures)  in all of these cases the sub-torus lies in a $\Gamma $ neighborhood. Therefore we have: 

\begin{Cor}
\end{Cor}
\begin{itemize}
\item[(a)] $S^{4}_{K}=S^4$\\
\item[(b)] $(S^{2}\times S^2)_{K}=S^2 \times S^2$\\
\item[(c)]  $[\;S^{1}\times S^3 \# (S^2 \times S^2)\;]_K =S^{1}\times S^3 \# (S^2 \times S^2)$
\end{itemize} 

By taking $K\subset S^3$ to be knots with different Alexander polynomials and using \cite{a1} we can state Corollary 2 (a) in the following useful form:

\begin{Thm} The fishtail F (the $2$-sphere with one self intersection) can imbed into $S^4$ infinitely many different ways $f_{K}: F \hookrightarrow S^4$, so that each $S^{4} -f_{K}(F)=F_{K}$ is a different  exotic copy of $F$, where $K$ are knots with different Alexander polynomials.
\end{Thm}

%\newpage

\section{ Proving $M(K)$ is standard}

\begin{Thm} $M(K)=S^{1}\times S^{3} \# (S^2\times S^2)$
\end{Thm}

\proof The first picture of Figure 8 is the handlebody of  $S^{1}\times S^{3}_{K}  $ surgered along the linking loop C (in the figure K is drawn as the trefoil knot), as discussed in  \cite{a2}. Here the pair of small red linking handles denotes the surgering the loop $C$ in $S^{1}\times S^{3}_{K} $. By sliding the zero framed circle over the $+1$ framed circle we obtain the second picture of Figure 8. Then by sliding the small  $-1$ framed red circle over one of the long zero framed circles (the ones going through the $1$-handle), and then sliding the large $-1$-framed circle over this small $-1$ framed circle we obtain the first picture of Figure 9 (now the large  $-1$ framed circle becomes $0$ framed). Note that this last move is from \cite{a2} (e.g. going from Figure 30 to Figure 31 of \cite{a2}). Then by  sliding $+1$ framed circle over the $-1$ framed circle we obtain the second picture of Figure 9 (i.e. the reverse of the first move of Figure 8). Now by \cite{a1}, this is just the handlebody of the knot surgered manifold of Figure 6, which is $[ \;S^{1}\times S^3 \# (S^2 \times S^2)\; ]_{K}$, so by the Corollary 2 it is $S^{1}\times S^3 \# (S^2 \times S^2)$. \qed.

 \begin{figure}[ht]  \begin{center}  
\includegraphics[width=.55\textwidth]{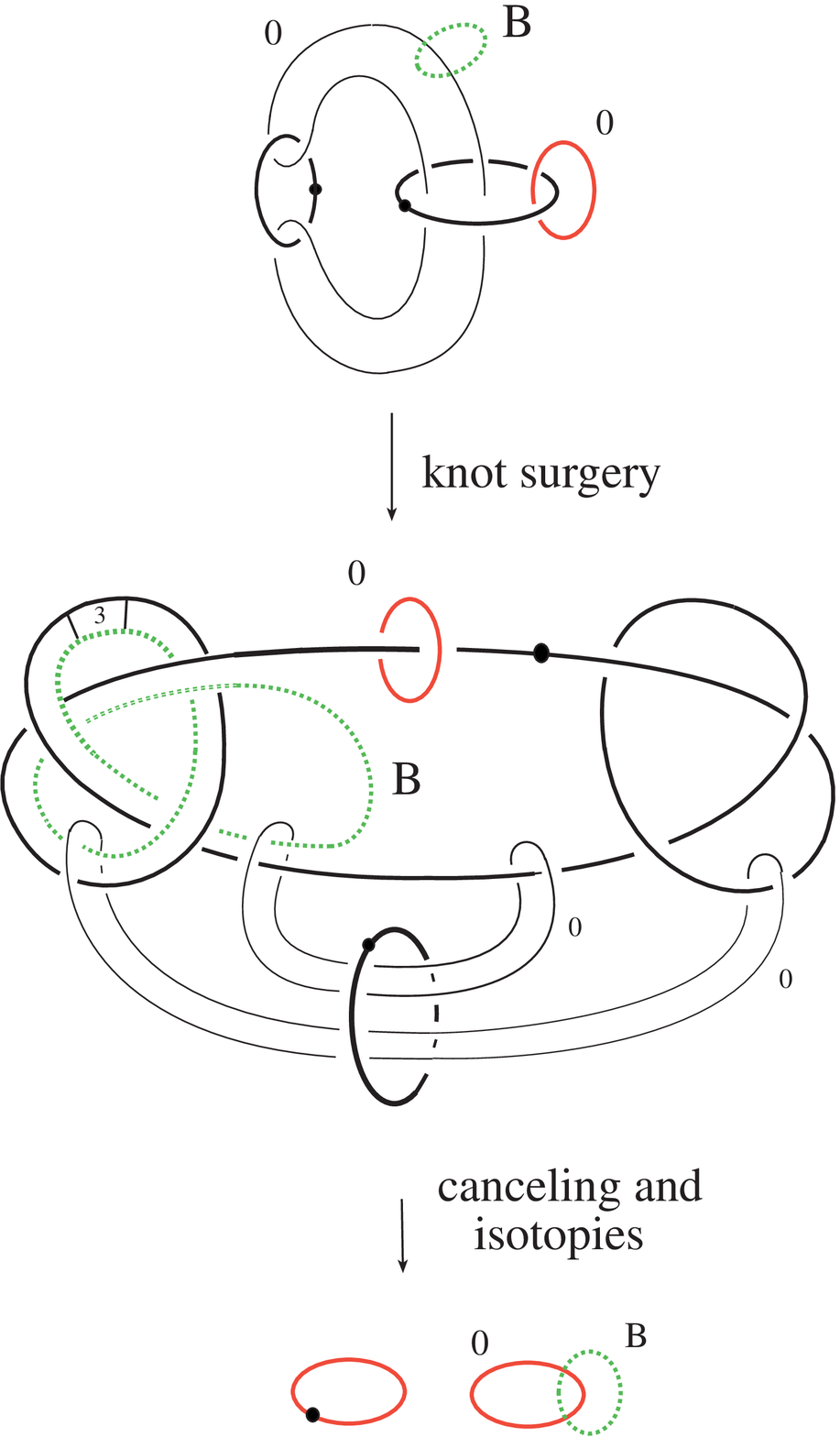}   
\caption{} 
\end{center}
\end{figure}

 \begin{figure}[ht]  \begin{center}  
\includegraphics[width=.6\textwidth]{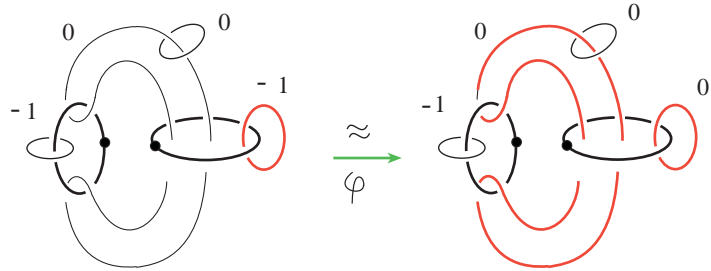}   
\caption{ $S^2 \times S^2$ as double of two cusps} 
\end{center}
\end{figure}

 \begin{figure}[ht]  \begin{center}  
\includegraphics[width=.9\textwidth]{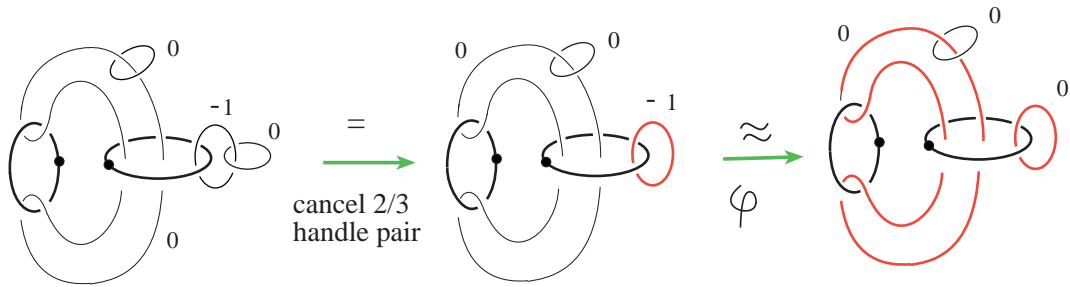}
\caption{$S^1\times S^3 \# (S^2\times S^2)$ as double of two fishtails}   
\end{center}
\end{figure}

 \begin{figure}[ht]  \begin{center}  
\includegraphics[width=.82\textwidth]{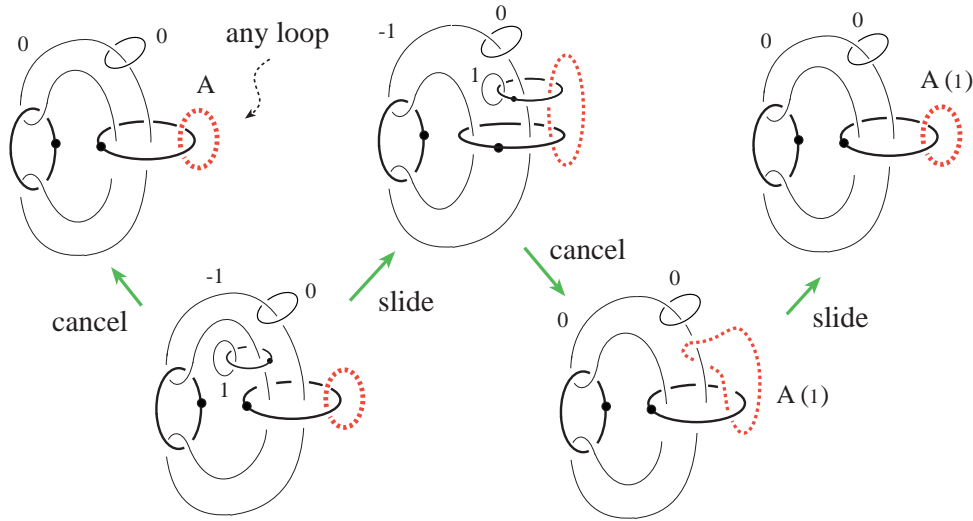}   
\caption{Diffeomorphism $\varphi: S^{2}\times T^2\to S^{2}\times T^2$} 
\end{center}
\end{figure}

 \begin{figure}[ht]  \begin{center}  
\includegraphics[width=.55\textwidth]{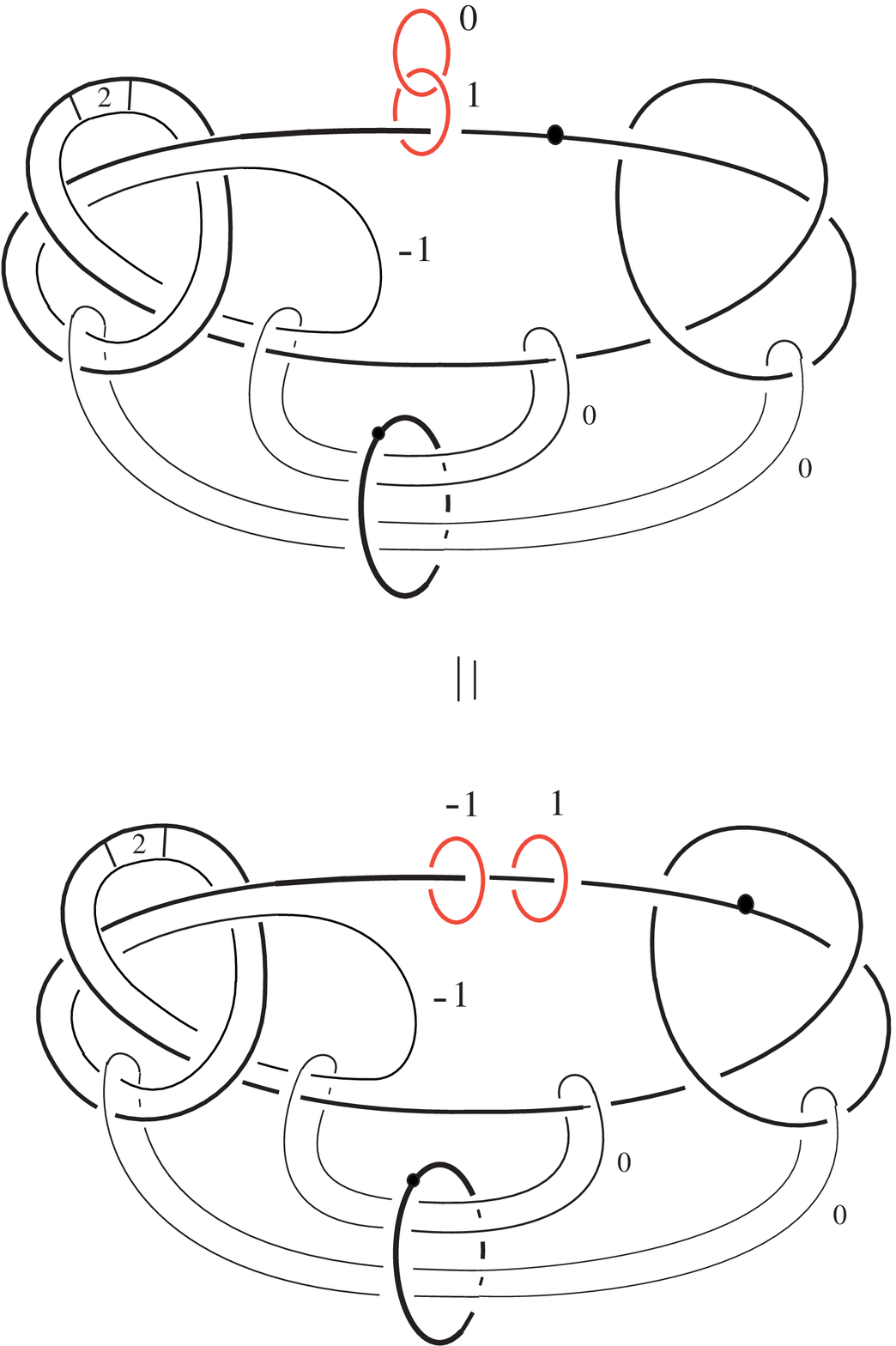}   
\caption{} 
\end{center}
\end{figure}

 \begin{figure}[ht]  \begin{center}  
\includegraphics[width=.55\textwidth]{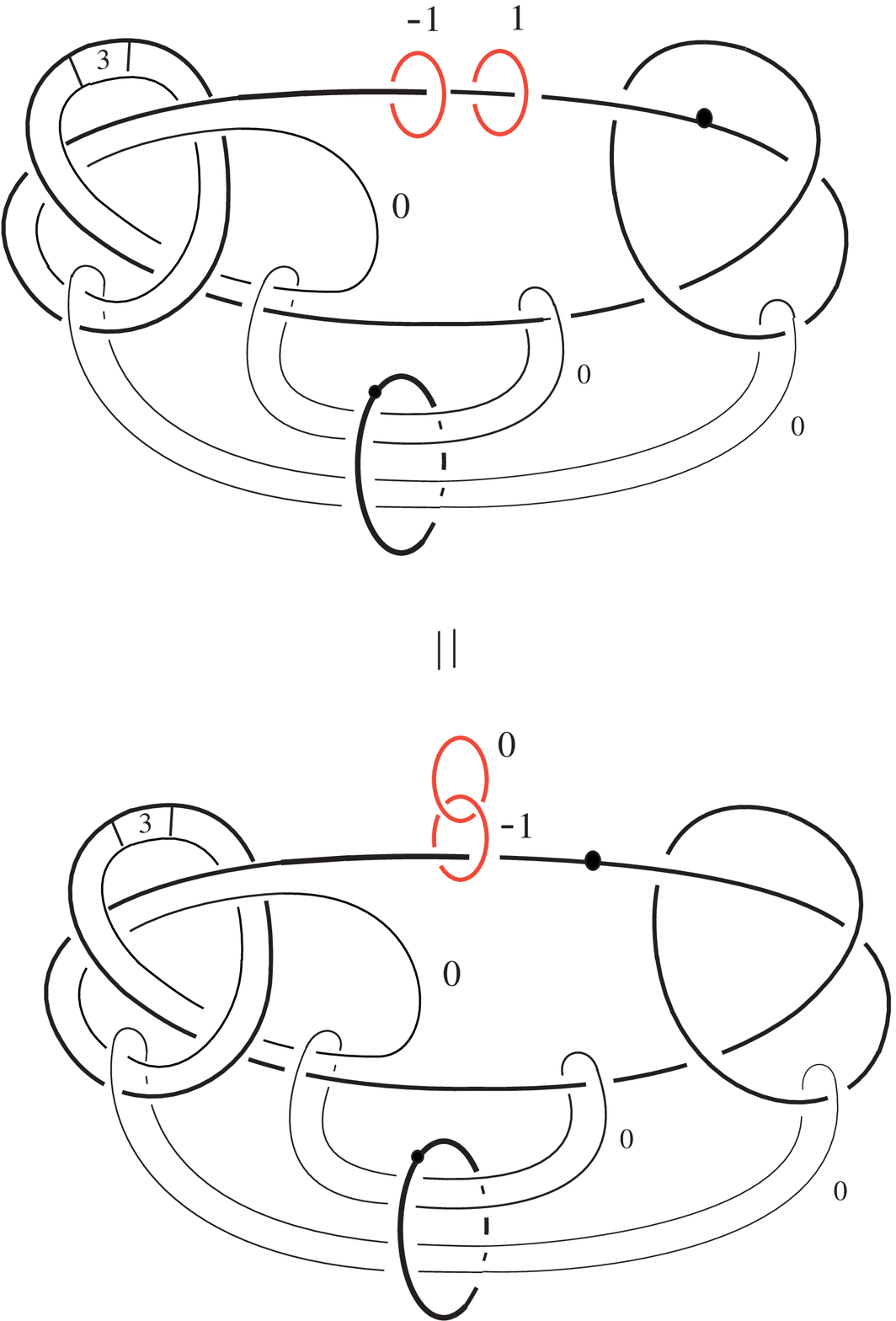}   
\caption{} 
\end{center}
\end{figure}

\end{document}